
\input harvmac
%

\def\({\left(}		
\def\){\right)}		
\def\[{\left[}		
\def\]{\right]}		
%
%
%

\def\la{\lambda}	
\def\de{\delta}

%

\lref\YB{D.Bernard, M.Gaudin, F.D.M.Haldane and V.Pasquier, J.Phys.A 26
(1993) 5219}
\lref\gaudin{M.Gaudin, La fonction d'onde de Bethe,
Edition du CEA, Masson 1983}
\lref\bour{N.Bourbaki, Groupes et algebres de Lie,ch 4,5,6}
\lref\mac{I.G. Macdonald, Seminaire N.Bourbaki, expose 797,1995}
\lref\cher{I. Cherednik, C.M.P 169, (1995) 441,461}
\lref\car{R.W.Carter Simple groups of Lie type (1972)}
\lref\lus{G.Lusztig, J.Amer.Math.Soc. 2 (1989), 599,635}
\lref\cours{V.Pasquier, Lecture notes in Physics 436 (1993), 36}

\Title{T95/093}
{{\vbox {
\bigskip
\centerline{Scattering matrices and}
\centerline{Affine Hecke Algebras.} }}}
\bigskip
\centerline{V. Pasquier}

\bigskip

\centerline{ \it Service de Physique Th\'eorique de Saclay
\footnote*{Laboratoire de la Direction des Sciences
de la Mati\`ere du Commissariat \`a l'Energie Atomique.},}
\centerline{ \it F-91191 Gif sur Yvette Cedex, France}

\vskip .5in

These lectures deal with two related topics.
The first one is the construction of the scattering matrix for
an arbitrary Weyl group.
The second topic is the application of the affine Hecke algebra
to construct physical models with enlarged symmetries.

\noindent
\vskip 2.0in
{\noindent \it Lecture at Schladming, March 4-11 1995 and at Trieste,
April 10-12 1995.}
\Date{07/95}

\newsec{introduction}
These lectures deal with two related topics.

The first one is the construction of the scattering matrix for
an arbitrary Weyl group. The aim is to define
a set of commuting operators in terms of generators which obey
the Yang-Baxter equation.
These operators which coincide with the usual scattering matrices
form a group isomorphic to the weight lattice of the corresponding
Lie group. In the first lecture, we  interpret the Yang Baxter generators
as elementary reflections. This geometrical picture
enables us to construct the scattering matrices in a systematic way
and to show that they commute.

In this lecture, we have used the language of Coxeter-Weyl groups
and cells exposed in \car, \bour.

The second topic is the application of the affine Hecke algebra
to construct physical models with enlarged symmetries.

The affine Hecke algebra is a deformation of the affine Weyl groups.
The reflections across the origin are deformed into the none affine generators
of the Hecke algebra and the translations are deformed into
generators which coincide with the scattering matrices
of the first lecture.
Using a specific representation of the Yang Baxter generators which
uses the Hecke algebra, we show that the scattering matrices
obey the affine Hecke relations.
We then derive a representation
in a space of polynomials in several variables.
This enables us to construct a
set of commuting Hamiltonians essentially given by the center of
the affine Hecke algebra. These Hamiltonians are Hermitian
and their spectrum can be computed in a simple way.
One of their main interests is that they commute with
q-deformations of the usual
affine Lie algebras defined in terms of
quadratic algebras.

The abstract theory of affine Hecke algebras is exposed in \lus.
We expand here the point of view of \YB, \cours.
Parallel arguments are in \cher
and some recent results can be found in \mac.

\newsec{The scattering matrix}


Two important tools of soluble models are
the transfer matrix an the scattering matrix.
In his study of the delta interacting gas, Yang introduced a commuting
family of scattering matrices. He obtained the spectrum
by diagonalizing simultaneously these matrices. A simple way to
prove the commutation relations for the
scattering matrices consisted in
defining a commuting family of transfer matrices depending of a continuous
parameter. By letting the parameter take some special values,
he recovered the scattering matrices. This approach is
expanded in the book of M.Gaudin \gaudin.

I do not know a generalization of the transfer matrix for an
arbitrary Weyl group, and when this can be done,
(A,B,C,D cases),
one cannot simply obtain the scattering matrix by letting the spectral
parameter
take special values except in the A case.
On the other hand,
the scattering matrix can be easily generalized to
arbitrary Weyl groups.
If a proper correspondence is used, the scattering matrices are
identified with translations in the affine Weyl group.
To generalize the construction of the scattering matrices
it is useful to
give them a geometrical interpretation.
It enables to
construct  commuting operators from a set of generators
obeying the (generalized) Yang-Baxter algebra.

In the first part, I
recall some known facts about Weyl groups following
the presentation of Carter \car.
I modify the
presentation of the Weyl group in terms of generators and relations
to introduce generators which
obey the Yang-Baxter equation.
The difference between the two approaches is analogous to the description
of a rigid body motion in the moving frame compared to the description
in the rest frame.

In the second part, I extend these considerations to the affine case
which enables to obtain the expression of the scattering matrix.
The geometrical interpretation is the following:
Periodic trajectories inside the fundamental cell associated to
a Weyl group are in correspondence with translations
in the weight lattice. To these trajectories, we associate
operators by taking the time ordered product of the
Yang-Baxter operators of
the reflections across the
walls hit by the trajectory. These operators
form a commutative group isomorphic to the weight lattice.


\subsec {Systems of roots and Weyl groups}
{\it Root systems}

Let $V$ be an euclidian space of dimension $l$. for $r$ a vector of $v$
we denote by $w_r$ the reflection in the hyperplane orthogonal to $r$.

\eqn\ref{w_r (x)=x-2{(r,x)\over (r,r)}r.  }

A system of roots is defined by a set $\Phi$ of non zero vectors
spanning $V$ which obey the following properties:

1)   if $r,s \in \Phi $ then, $w_r(s) \in \Phi$.

2)   if $r,s \in \Phi $ then  $2(r,s)/(r,r)$ is a rational integer.

3)   if $r,\la r \in \Phi$ then $\la=\pm 1$

It will be convenient to denote $r^{V}=2r/{(r,r)}$.
The lattice spanned by the $r^{V}$ is denoted $Q^{V}$.

We shall also consider the lattice spanned by the vectors $p$ such that
$ ( p,r)$ is  a rational integer for any root $r\in \Phi$.
This lattice denoted $P^{V}$ contains the coroot
lattice $Q^{V}$.

One can define a basis of $V$ , $\Pi$, which is a subset of $\Phi$
and is such that any root of $\Phi$ is  a linear combination of roots of $\Pi$
with coefficients which are either all non-negative or all non-positive.
$\Pi$ is called a fundamental system of roots.
$ \Pi = \{r_1,r_2,...,r_l\} $
and a root can be written:
\eqn\pos{ r=\sum_{i=1}^l \la_i r_i }
where either $\la_i \ge 0$ for all $i$, or $\la_i \le 0$ for all $i$.
Accordingly, we say that $r$ is a positive $r\in\Phi^+$ or a negative root
$r\in\Phi^-$.
We say that $r>s$ if $r-s$ is a linear combination with coefficients
$\ge0$ of elements in $\Pi$.
Every positive systems
of roots contains just one fundamental system. Thus there is
a one to one correspondence between
positive systems and fundamental systems in $\Phi$.

{\it chambers}

For our purpose, it is useful to have a geometric interpretation of
a fundamental system.
For each root $r\in\Phi$, we denote by $H_r$ the hyperplane orthogonal to $r$.
The set
$V-\cup H_r$
is disconnected, its connected components are called chambers.
The roots orthogonal to the bounding hyperplanes of a chamber and pointing
into the chamber form a fundamental system $\Pi$.
A chamber $C_{\Pi}$ is thus defined as the
set of points $x\in V$ such that:
\eqn\cha{ (r_i,x)>0, \forall r_i\in \Pi }
Moreover, every fundamental system
arises in this way from some chamber.

Let $C$ be a chamber, and $\de(C)$ be its boundary. The bounding hyperplanes of
$C$ are called the walls of $C$ and their intersection with $\delta(C)$ the
faces of $C$. Two chambers which have a face in common are called adjacent.

{\it Weyl group}

The group $W$ generated by the reflections $w_r$ is
called the Weyl group of $\Phi$.
Each element of $W$ transforms $\Phi$ into itself.
One can show that the Weyl group is generated by the
fundamental reflections $w_r$ with $r \in \Pi$.
Given two chambers $C$ and $C'$ (alternatively two fundamental systems
$\Pi$ and $\Pi'$), there is a unique element of the Weyl group
such that $w(C)=C'$.
Let us here give a construction of this element.
This presentation will make the  Yang-Baxter equation appear naturally
in this context.

We consider a sequence $C_0,C_1,...,C_m$ of chambers
with $C_0=C,\ C_m=C'$ such that two
consecutive chambers are adjacent. Each chamber is characterized by a
fundamental system $\Pi_i$, and there is a unique $r_i \in \Pi_i$ such that
$w_{r_i}C_i=C_{i+1}$.
We define the element of the Weyl group such that $w(C)=C'$ by:
\eqn\douw
{ w=w_{r_{m-1}}w_{r_{m-2}}...w_{r_0} }

{\it length of the Weyl group elements}

The length of the element $w$, $l(w)$ can be defined in two different ways.
The first one being:

a) $l(w)$ is the number of positive roots turned by $w$ into negative roots.

Equivalently, if we denote $\Phi^+,\ \Phi'^+$ the positive roots
respectively associated to the fundamental systems of $C$ and $C'$.:

\eqn\len{ l(w)=|\Phi^+ \cap \Phi'^-|  }
Let $r\in\Pi'$. Consider the word $w_rw$, since $w_r$ changes the sign
of only $\pm r$ in $\Phi'$,
one has in particular:
\eqn\lenn{
\eqalign{
l(w_rw)&=l(w)+1 {\rm\ \  if\  \ } r\in \Phi^+ \cr
l(w_rw)&=l(w)-1 {\rm\ \  if\  \ } r\in \Phi^-  \cr
}}
Now, consider an expression of the form \douw~
to represent the Weyl group element $w$. From \lenn~ it follows that
the minimal number of reflections in \douw~ is
larger than l(w). Let us assume that the number of terms is strictly
larger than l(w). This means that the expression contains a
reflection $w_{-r}$ with $r\in\Phi^+$. From the above argument, this
is only possible if the reflection $w_r$ has occured earlier.
Thus, the expression is of the form:
\eqn\inu{
w=w_{r_{m-1}}...w_{-r}Xw_r...w_{r_0} }
Let us denote $w_r(X)$ the expression obtained from $X$ by substituting
$w_{w_r(p)}$ for $w_p$ in the expression of X. A shorter expression for
$w$ is given by:
\eqn\inut{w=w_{r_{m-1}}...w_{r}(X)...w_{r_0}}
It follows that another definition of the length is:

b) l(w) is the the minimal length of an expression $w$ in the form $\douw$.

A word of minimal length is one for which all the $r_i$ in the
expression \douw~of $w$ are in
$\Phi^+ \cap \Phi'^-$.

Given two chambers $C$ and $C'$, let us describe the way to construct
an element $w$ of minimal length such that $w(C)=C'$.
For this, we proceed by induction on $l$. If $l=0$, $C=C'$ and $w=1$.
If $l\ne 0$, $\Phi'^-\cap \Pi$ is not empty. Let $r_0$ belong to this set.
We put $w_{r_0}(C)=C_1$ and $\Pi_{1}=r_0{\Pi}$. The set of positive roots
defined
by $\Pi_{1}$ is obtained from $\Phi$ by replacing $r_0$ by $-r_0$.
Hence
\eqn\long{|\Phi^+_1 \cap \Phi'^-|=l-1 }
and we can continue the construction replacing $C$ by $C_1$.

{\it Yang-Baxter equation}

We now give a presentation of the Weyl group
$W$ as an abstract group in terms of generators $x_r$, $r\in \Phi$
using the above construction.
Namely, we consider the elements \douw~with the reflection $w_r$
replaced by $x_r$ and we impose relations to the $x_r$ in such a way
that two elements $w_1$ and $w_2$ coincide whenever the Weyl group
reflections are equal.

The operators $x_r$ are subject to the following relations:

a) unitarity:
\eqn\yba{x_{-r}x_r=1}

b) Generalized Yang-Baxter equation:

Let $r$ and $s$ be two roots in some fundamental system $\Pi$,
and $m_{rs}$ be the order of $w_rw_s$. Then, we have:
\eqn\ybb{
\eqalign{&x_{r_1}x_{r_2}x_{r_3}...x_{r_m}=
x_{r_m}x_{r_{m-1}}x_{r_{m-3}}...x_{r_1}\cr
{\rm where:\ }
&r_1=s,\ r_2=w_s(r),\ r_3=w_sw_r(s),...,\ r_m=r \cr}}

For example, in the case where the root system is $A_n$, if we take
$r=e_1-e_2$,
$s=e_2-e_3$, then $m_{rs}=3$ and the above equation writes:
\eqn\usu{x_{e_2-e_3}x_{e_1-e_3}x_{e_1-e_2}=
x_{e_1-e_2}x_{e_1-e_3}x_{e_2-e_3} }
one recognizes the usual Yang-Baxter equation.
These equation are generalization of the Yang-Baxter equation for arbitrary
Weyl groups.

\subsec {Affine root systems and Weyl groups}

{\it affine-root systems}

In order to define a scattering matrix, we need
to extend the previous considerations to
the affine case.

We denote by $E_0$ an affine vector space
of dimension $d$.
For $r\in \Phi$ and $k\in Z$ we define a hyperplan in $E_0$ by:
\eqn\plan{H_{r,k}=\{x \in E_0\vert (r,x)=k\}}
The orthogonal reflections with respect to the Hyperplanes $H_{r,k}$
generate a group called the affine Weyl group.
This group is also the semi direct product of the Weyl group defined
earlier and the translations in the coroot lattice $Q^{V}$.

The affine-roots are defined by a couple of non affine root and a rational
integer. We denote $\tilde\Phi$ the set of affine roots.
Each affine root $r+k\de$
defines the Hyperplane $H_{r,k}$ and the corresponding
reflection in the affine Weyl group:
\eqn\aref{w_{r+k\de}(x)= x-((r,x)-k)r^{V}}
A fundamental system $\tilde\Pi$ is given by a basis of $\tilde\Phi$ such that
any root of $\tilde\Pi$ is a linear combination of roots of $\tilde\Pi$
with coefficients which are either
all positive $(\tilde\Phi^+)$  or all negative $(\tilde\Phi^-)$.
For example, given a fundamental system $\Pi$,
if we denote $r_{m}$ is the largest root of $\Phi$,
we obtain a
fundamental system $\tilde\Pi$ given by:
\eqn\fonda{
\tilde\Pi=\{r|r\in \Pi\} \cup \{-r_{m}+\de\} }

The positive roots are given by:
\eqn\posi{\tilde\Phi^+= \{r+k\de|r\in\Phi,k>0\}\cup\{r|r\in\Phi^+\}}

{\it Cells}

The set $E_0-\cup H_{r,k}$ is disconnected, its connected
components are called cells.
A cell $C_{\tilde\Pi}$,
$\tilde\Pi=\{r_i+k_i\de,\ i=0,...,r\}$,
is characterized by the fundamental system  given by the roots
which define its bounding hyperplanes and which point into it:
\eqn\poedt
{C_{\tilde\Pi}=\{x\in E_0\ |(r_i,x)\ge k_i,\ \forall i\} }

The faces of a cell are defined as earler,
and two cells which have a face in common are called adjacent.

{\it representation of the Weyl group elements by words}

We denote by $E_p$ the affine vector space centered at some point $p$.
$E_0=\{0,x\}$ , $E_p=\{p,x\}$, with the identification $\{0,x\}=\{p,x-p\}$.
To simplify the notations, $\{0,x\}$ is denoted $x$.
In the following, we shall restrict to the affine spaces $E_p$
with $p$ in the coweight lattice: $p\in P^V$.

Let us define operators which intertwine $E_p$ and $E_{p+a}$
and which are equal to one when we use the equivalence relation
to identify the two spaces. They simply translate the origin of the
affine space.
They are defined by:
$$t^a\{p,x\}=\{p+a,x-a\}$$
We also consider the affine reflections $x_r$ for $r\in\Phi$ acting in $E_p$
as:
$$x_r\{p,x\}=\{p,w_r(x)\}$$
Inside the parentheses, $w_r$ denotes the Weyl reflection defined earlier.
Using the equivalence relation to identify $E_p$ with $E_0$,
the reflection
$x_r|_{E_p}
$ is identified with $w_{r+(r,p)\de}$.

Given two cells, C and C', there is a unique element $w$
of the affine Weyl group
which maps the first one into the second one.
Let us here construct a word representing this element using
the elementary transformations $x_r$ and $t^a$.
We consider a sequence of cells
$C_0=C,C_1,... ,C_m=C'$ such that two consecutive cells
are adjacent. To simplify the discussion, we fix
$C$ to be defined by the fundamental system \fonda.
We consider a linear transformation acting in $E_0$ as follow:
\eqn\word{
S=t^{a_{m}}x_{r_{m-1}}t^{a_{m-1}}x_{r_{m-2}}...t^{a_1}x_{r_0}t^{a_0} }
The reflection $x_{r_i}$ acts in $E_{a_{i}+a_{i-1}+...}$ and maps
the cell $C_i$
onto the cell $C_{i+1}$.
The origin $a_{i}+a_{i-1}+...$ must therefore belong
to the common wall of the cell $C_i$ and $C_{i+1}$.
The root $r_i$ is orthogonal to this wall
and we take it to point into the cell $C_i$.
It is clear that $S$ intertwines
$E_0$ with $E_p$ where:
$p=a_0+a_1+...+a_m$.

Using the
equivalence relation to identify $E_0$ with $E_p$,
one can view the
transformation $S$ in the affine Weyl group.
Conversely any element
$w$ in the affine Weyl group
is
characterized by a cell $C'$ and
can be represented by a word  \word.
We can always choose $p\in Q^V\cap C'$, in this way
$w$ characterizes a unique point $p\in Q^V$ which is a summit of $C'$.
In this case, we can identify the Weyl group element $w$ with
the word denoted $S_w$ which represents it.

Here, we need to be more general and we take
$p\in P^V\cap C'$. Unless $P^V=Q^V$, this does not
determine a unique point $p$.
In that case, we denote by $S_{w,p}$ the word \word~
which intertwines $E_0$ and $E_p$ and which coincides with the element
$w$ in the affine Weyl group.
$p$ is a summit of $C'$ and can be
written in a unique way as
\eqn\defgam {p=q-\gamma }
where
$q\in Q^V\cap C'$.

{\it Scattering matrix}

Let us restrict to the
case where $C'$ is obtained from
$C$ by a translation :
\eqn\defcell{C'=C+\xi   }
In order for $C'$ to define a cell,
$  \xi$ must be in $P^V$.
We denote by $w_{\xi}$ the corresponding element in the affine Weyl group.
($w_{\xi}$ is not a translation unless $\xi$ is in the coroot lattice
$Q^{V}$).
In the construction of $S_{w_{\xi},p}$,
one can choose $p=\xi$.
We denote by
$S_{\xi}$ the corresponding word $S_{\xi}=S_{w_\xi,\xi}$

In this case, $\gamma$ defined in \defgam~must be a summit of $C$;
therefore,  either $\gamma=0$, or $\gamma\in P^V$ and $(\gamma,r_m)=1$
where $r_m$ is the largest root of $\Phi$.
Such weights are called minuscule weights and are representative of
$P^V/{Q^V}$
in the weight lattice.

Let us now consider that the sequence of cells $C_i$ are
defined up to a translation. We
can view $S_{\xi}$ as intertwinning $E_p$ and $E_{p+\xi}$
for $p\in P^V$  arbitrary .
Then, it is clear that two words $S_{\xi}$ and $S_{\xi'}$  can be
multiplied and one has:
\eqn\mulyip{S_{\xi_1}S_{\xi_2}= S_{\xi_1+\xi_2}   }
Thus these linear transformations form a commutative
group isomorphic to the translations in the coweight lattice.
Note that the above relation is not true for the
Weyl group elements:
$w_{\xi_1}w_{\xi_2}\ne w_{\xi_1+\xi_2}$
unless $\xi_2$ is in the root lattice.

{\it remark}

We can be slightly more general and  define
a semi group multiplication law on the words $S_{C_1,p_1}^{C_2,p_2}$
which transform $C_1,p_1$ into $C_2,p_2$ by defining the product:
\eqn\produi{
S_{C_2,p_2}^{C_3,p_3}
S_{C_1,p_1}^{C_2,p_2}=
S_{C_1,p_1}^{C_3,p_3} }

{\it length}

We can define the length of
a Weyl group element $w$ as the minimal length of $S_{w}$
in terms of operators $x_r$. The length is also given by the expression
which generalizes \len,\lenn~ in an obvious way.
\eqn\lena{ l(w)=|\tilde\Phi^+ \cap \tilde\Phi'^-|  }
where $\tilde\Phi'^-=w(\tilde\Phi^-)$.
The same definition applies to the generalized words $S_{w,p}$
since $p$ is not relevant in the definition of the length.

In the case where the  word \word~is of minimal length, one has:
\eqn\noulong{
\tilde\Phi^+ \cap \tilde\Phi'^- =
\{r_i+\delta\sum_{j=0}^{i-1}(r_i,a_j)\} }

Consider the dominant weights
defined by:
\eqn\dom{P^{V}_{\rm dom}=\{\xi\in P^{V}|(\xi,r)\ge0\ ,\forall r\in\Pi\}}
Let us restrict to the words $S_{\xi}$ where $\xi$ is a dominant weight:
$\xi\in P^{V}_{\rm dom}$.
In this case, it is not difficult to see that:
\eqn\nougiu{
\tilde\Phi^+ \cap \tilde\Phi'^- =
\{r+\delta k|\ \  r\in \Phi^+,\ \ 0\le k<(\xi,r)\} }
and we have:
\eqn\geur{l(w_{\xi})=\sum_{r\in\Phi^+}(\xi,r)}
Moreover, it follows from \noulong~that
the reduced expression of $S_{\xi}$, the oprator $x_{r}$ occurs
exactly $(r,\xi)$ times and the operator $x_{-r}$ never occurs.

It is also straightforward to construct a word of minimal length by extending
the method used in the non affine case.
A geometrical interpretation of this construction is the following:
One draws a straight line joining an arbitrary point $A$
in the cell $C$ to its translated $A'=A+\xi$ in $C'$.
Traveling along this geodesic one goes successively  through the cells $C_i$.
The word $S_{\xi}$ is obtained taking the product (from the right to the left)
of the reflections across the hyperplanes which are successively
come through
by the trajectory.
The translations $t^a$ must be inserted in order for the origin of the affine
space to belong to the reflecting hyperplane and one must have
$\sum_{i=0}^m a_i=\xi$.

{\it definition of words by generators and relations}

We now proceed as before and replace  the reflection operators
$x_r$ and the translation operators $t^a$ by abstract operators.
We impose  relations to these operators  so
that the two words $S_{w,p}$ and $S_{w',p'}$ are equal when:

a) They both intertwine $E_0$ with the same affine space $E_p$.

b) After identifying $E_0$ with $E_p$, the two affine Weyl
group elements $w$ and $w'$ coincide.

For this we impose the conditions \yba~ and \ybb~ to the operators $x_r$
and for the operators $t^a$ we require:

1) that they form an additive group:
\eqn\yta{t^at^b=t^{a+b} }

2) that $t^a$ commutes with $x_r$ whenever $r$ and $a$ are orthogonal:
\eqn\ytb{t^ax_r=x_rt^a,\ {\rm if}\ (a,r)=0 }

The elements $S_{\xi}$ constructed in that way clearly commute, we
call them scattering matrices.
Using this construction we have obtained a method to construct a group of
commuting operators isomorphic to the coweight lattice.

{\it example:}

Let us  for example consider the case where the Weyl group
is $A_n$. If we take the vector $\xi=-e_i$,
the scattering matrix is given by:
\eqn\exaa{S_{-e_i}=x_{e_i-e_{i-1}}x_{e_i-e_{i-2}}...x_{e_i-e_1}t^{-e_i}
x_{e_i-e_n}....x_{e_i-e_{i+1}}  }
In the case where the $t^a$ are equal to $1$,
these coincide with the scattering matrices considered by Yang.


In the next section, we shall use these transfer matrices to construct
the q analogue of the Dunkl operators and we shall show that they obey the
defining relations of a affine Hecke algebra.

\newsec{The affine Hecke Algebra}
In this section, we use the previous formalism to construct some
representations of the affine Hecke algebra.
In the first part, we recall some well known results about
representations of the Yang-Baxter algebra using generators
which obey the (non affine) Hecke algebra relations. This gives a certain
representation of the scattering matrices.
In the second part, we consider a limiting form of these scattering matrices
to obtain a representation of the affine Hecke algebra.
In the third part, we derive a representation of
the Hecke and of the affine Hecke
algebra acting in the group
algebra of the weight lattice.
We  obtain
in this way the q-Dunkl operators.

\subsec{representation of the Yang-Baxter operators}

{\it Hecke algebra}

Let us consider a set of operators $g_{r}$ indexed by the roots $r\in\Phi$
on which the Weyl group acts in a natural way:
\eqn\cov{w_rg_s=g_{w_r(s)}w_r}
and
which obey the following relations:

a) Braid group relations:

Let $r$ and $s$ be two fundamental roots in $\Pi$, and $m_{r,s}$ be
the order of $w_rw_s$, then:
\eqn\braid{g_rg_s... =g_sg_r...
,\ m_{r,s}\ {\rm terms\  on \ either \ side}  }

b) Hecke relations

\eqn\heca {(g_r-q_r)(g_r+q_r^{-1})=0  }
where $q_r$ is a complex number which depends on length of the root $r$.

The subset $g_r$ with $r\in\Pi$ generates the Hecke algebra associated to
$\Phi$ defined by:
\eqn\hecb{g_wg_{w'}=g_{ww'} \ \ {\rm if\ }\ l(ww')=l(w)+l(w')  }
Note that we have $g_{w_r}=g_r$.
only for  $r\in \Pi$.

{\it Group algebra of the weight lattice}

The the group
algebra of the Weight lattice is denoted $P$
and is called the space of spectral parameters:
\eqn\wela{
P=\{\sum c_{\la}e^{\la} | \la \in P^V\} }
with the product given by:
\eqn\pols{
e^{\mu}e^{\la}=e^{\la+\mu} }
The  Weyl group is not acting on $P$: $w_r e^{\la}=e^{\la} w_r$.

{\it Yang-Baxter generators}

It is straightforward to verify that the following operators
obey the unitarity and the Yang-Baxter equation \yba,\ybb:
\eqn\rey{
x_r=w_r{e^{r}g_r-g_r^{-1}\over
e^{r}q_r-q_r^{-1}}  }


{\it relations satisfied by the} $t^a$:

In order for the operators $t^a$
to obey the relations \ytb , we require
that they commute with $P$ and satisfy the following relations:
\eqn\nyta{\eqalign{
wt^a&=t^{w(a)}w,\ \forall w\in W\cr
t^ag_r&=g_rt^a,\ {\rm if}\ (r,a)=0\cr }}

{\it Remarks}:

a)
Another way to realize \ytb~ consists in having
the operators $t^a$ commute with the Weyl group
and obey the following relation with $P$:
\eqn\spii{
t^a e^r=e^{(r,a)} e^r t^a }
such a realization occurs in the so called quantum Kniznik Zamolodchikov
equation.

b)
Let us consider the case of a root system
for which the third condition defining a root system
is not satisfied ( nonreduced system ) and there
are proportional roots $r$ and $2r$.
Those are the BC systems for which
the roots are: $\pm e_i,\ \pm 2e_i,\ \pm e_i\pm e_j$.
In this case. one can find an expression for $x_r$ which interpolates between
the two forms of \rey~for $x_r$ and $x_{2r}$.
\eqn\reyb{
x_r=w_r{e^{2r}g_r-g_r^{-1}+e^{r}\alpha\over
e^{2r}q_r-q_r^{-1}+e^{r}\alpha }  }
This formula coincides with $x_{2r}$ in \rey~
if $\alpha=0$ and  $q_{2r}=q_r$. It coincides with $x_{r}$ if
$\alpha=q_r-q_r^{-1}$. It also satisfies the unitarity condition \yta~ and
the Yang-Baxter equation \ytb.


\subsec{affine Hecke relations}
In this section, we consider the Hecke algebra \hecb.
We complete it by adding a group of translations isomorphic to
the coweight lattice generated by the $S_{\xi}$.
Moreover, we require that the translations and the generators $g_r$
obey the commutation relations which define a affine-Hecke algebra.
In addition to the relations \hecb~ satisfied by the $g_r$, $r\in\Pi$
and the  additive group formed by the $S_{\xi}$, $\xi \in P^V$, we require that
the $g_r$ and the $S_{\xi}$ obey the relations:

Let $r\in\Pi$, and $\xi \in P^{V}$, then:
\eqn\heka{
\eqalign{
i)\ &{\rm if}\ (r,\xi)=0,\ g_rS_{\xi}=S_{\xi}g_r \cr
ii)\ &{\rm if}\ (r,\xi)=1,\ g_rS_{\xi}=S_{w_r(\xi)}g_r^{-1} \cr  }}

Here we consider a limit of the operators $S_{\xi}$ constructed in \word~
which obey these relations.
We consider the following limiting form (
obtained when we set $e^r=\infty $ for $r>0$
with respect a fixed fundamental system $\Pi$)
of the operators
$x_r$ \rey~:
\eqn\ixr{
\eqalign{ x_r&=w_rg_r \ \ {\rm if}\ r\in\Phi^+ \cr
x_r&=w_rg_r^{-1} \ {\rm if}\ r\in\Phi^- \cr }}

Let us first show these relations in the case where $\xi$
is a dominant weight, $\xi\in P^{V}_{dom}$.
We recall that in this case, $S_{\xi}$ has a reduced expression:
\eqn\wordxi{
S_{\xi}=t^{a_{m}}x_{r_{m-1}}t^{a_{m-1}}x_{r_{m-2}}...t^{a_1}x_{r_0}t^{a_0} }
where $r_i\in \Phi^+$ and $x_{r_i}$ occurs $(r_i,\xi)$ times
in this expression.

Given a word $S$ and $w\in W$, we denote by $w(S)$ the word obtained by
substituting everywhere $t^{w(a)}$ and $x_{w(r)}$ for $t^a$ and $x_r$
in the expression of $S$.
In general, $w(S)\ne wSw^{-1}$ because there can be some $x_{r_i}$ in
the expression of $S$ for which $r_i\in\Phi^+$ and
$w(r_i)\in\Phi^-$. The main point of the following proof is to
bring us back to a situation where the equality applies.

$i)$ {\it Proof of the first relation} \heka:

We consider the two words $x_rS_{\xi}$ and $w_r(S_{\xi})x_r$.
Both words intertwine $E_0$ with $E_{\xi}$ and after identifying the
two spaces, they coincide with the Weyl transformation
$w$ defined by its action on the cell $C$ defined by \fonda:
$w(C)=w_r(C)+\xi=w_r(C+\xi)$. Thus they are equal and one has:
\eqn\inter{w_rg_rS_{\xi}=w_r(S_{\xi})w_rg_r }
Now, since $(\xi,r)=0$ and $S_{\xi}$ is reduced,
$S_\xi$ contains only operators $x_{r_i}$ with
$r_i \ne r$ and $w_r(r_i)$ are all in $\Phi^+$.
Thus,
$w_r(S_{\xi})=w_rS_{\xi}w_r$. Substituting this equality in \inter,
$i)$ follows.

$ii)$ {\it Proof of the second relation} \heka:

Since $(\xi,r)=1$, one can write $S_{\xi}$ in the form:
$S_{\xi}=S'x_{r}$. $S'$ is a word of length
$l(S')=l(S_\xi)-1$ and $x_r$ does not appear in the expression of $S'$.
Using the geometrical interpretation, one shows that:
\eqn\intera{x_{-r}w_r(S')=S_{w_r(\xi)} }
Substituting $w_r(S')=w_rS'w_r ,\ x_{-r}=g_r^{-1}w_r$ and $S'=Sx_{-r}$ in
\intera~ $ii)$ follows.

If $\xi$ is not a dominant weight, we can always write $\xi=\xi_1-\xi_2$
where $\xi_1$ and $\xi_2$ are dominant weights such that:
$(\xi_1,r)=0$ and $(\xi_2,r)=0$ in case $i)$;
$(\xi_1,r)=1$ and $(\xi_2,r)=0$ in case $ii)$.
Writing $S_{\xi}=S_{\xi_1}S_{\xi_2}^{-1}$ we obtain the general result.

\subsec{Another presentation of the affine-Hecke algebra}

Consider the affine generator  given by:
\eqn\gross{g_{-r_{m}+\de}=t^{-r_m/2}g_{-r_m}t^{r_m/2} }
We could have considered the algebra
generated by the subset $g_r$ with $r\in \tilde \Pi$
defined by the relations \braid~\heca , replacing the non affine root system
by the affine root system $\tilde \Pi$.
Commuting the Weyl group elements $w_r$
to the left
of the word $S_{\tilde w}$ and reorganizing the translations $t^a$  so as
to make only the generators $g_r$ with $r\in \tilde \Pi$ appear to the
right
we obtain the following expression:
\eqn\zut{
S_{\tilde w}= w^{-1}g_{\tilde w} }
where  $g_{\tilde w}$ is in the subalgebra
generated by the $g_r$,$r\in\tilde \Pi$ and
$w$ is the projection of $\tilde w$ in the Weyl group.

In the general case,
we obtain the following expression:
\eqn\zutu{
S_{\tilde w,p}= t^{-\gamma}w^{-1}g_{\tilde w}}
$\gamma$ is defined in \defgam.
In particular:
\eqn\zuta{
S_{\xi}= t^{-\gamma}w_{\xi}^{-1}g_{\tilde w_{\xi}}  }
in this case $\gamma$ is a minuscule weight.
The left-hand side operator
$ t^{-\gamma}w_{\xi}^{-1} $ depends only of $\gamma$, the projection of
$\xi$ in $P^V/Q^V$. These operators can be added to the generators
$g_r$ with $r\in \tilde \Pi$ to give a presentation of the affine
Hecke algebra by generators and relations.

\subsec{representation of the Hecke-algebra in a polynomial space}

In this section, we construct a representation of the scattering matrix
\word~
in a space of polynomials in several variables.
If we use the expression \rey~ for the operators $x_r$,
all we need is a representation of the operators $g_r,\ r\in \Pi$,
such that the Hecke-algebra relations
\cov, \braid, \heca~are satisfyed.
We also need a representation of the operators $t^a$
satisfying \nyta.
Such a representation of the Hecke algebra
is known in the mathematic literature (Lusztig,
Lascou and Schutzenberger). Our aim here is to deduce it from
the Yang-Baxter equation. In doing so, we shall discover
two different sets of operators which obey the affine-Hecke relations:
One is given by the $S_{\xi}$, the other by the polynomials on which the
$g_r$ act.

To construct the representation of the $g_r$,
We recall that the the group
algebra of the Weight lattice is denoted $P$.
There is a natural action
of the Weyl group on $P$  given by:
\eqn\wepol{
se^{\la}=e^{s(\la)}s }

We consider
the Hecke algebra \hecb~generated
by $f^r,\ r\in\Pi$.
Here the $f^r$ commute with
the Weyl group action and the multiplication by $e^{\la}$.

We define the operators $y_r$ as follows:
\eqn\reyy{
y_r=s_r{e^r f_r-f_r^{-1}\over
e^r q_r-q_r^{-1}}}
It is easy to show that the Yang-Baxter equations \yba~\ybb~
are equivalent to the fact that
$y_r$ with $r\in \Pi$ obey the defining relations of the generators
of the Weyl group:
\eqn\uni{y_r^2=1}
\eqn\unii{(y_ry_s)^{m_{rs}}=1}
where $m_{rs}$ is the order of $w_rw_s$.
It suffices for that to write $x_r=w_rz_r$ in \rey~ and
to  commute the Weyl group elements $w_r$  to the
left of the expressions \yba,\ybb. If one does the
same by commuting $s_r$ to the left of \uni,\unii, the
identities to verify are the same in both cases.

Let us now quotient the group algebra: $\{\sum c_{\la,w,w'}e^{\la}s_wf_{w'}\}$
by the relation $y_r=1$ to the right.
We denote by $\pi(.)$ the projection
which consists in eliminating the reflections $s_r$ to the right of an
expression in the quotient.
This operation eliminates the $s_r$ in the following way:
\eqn\subs{
\pi(...s_r)=\pi(...){e^r f_r-f_r^{-1}\over
e^r q_r-q_r^{-1}} }
The consistency of this projection
is assured by the relations \uni,\unii.
Alternatively, we can define the projection $\pi$ as follows:
\eqn\repga{
\pi(...(q_r s_r+(q_r-q_r^{-1}){1\over e^r-1}(s_r-1)))
=\pi(...)f_r }
Let us rewrite the above relation as
$\pi(...g_r)=\pi(...)f_r$.
The consistency of this operation is assured by the fact that the $g_r$
obey the Hecke algebra relations.
The expression for $g_r$ is
given by:
\eqn\repg{
g_r=q_r s_r+(q_r-q_r^{-1}){1\over e^r-1}(s_r-1) }
and the expression \ixr
of $x_r$ is:
\eqn\repx{
x_r=q_r +(q_r-q_r^{-1}){1\over e^{-r}-1}(1-s_r) }

One can verify that this representation of the
$g_r$ for $r\in\Pi$ and the translations $e^{\la}$
obey the defining relations \heka~ of a affine Hecke algebra.
In fact, it is easy to see that the relations \heka~ are equivalent
to the following relations:
\eqn\repgg{
g_rQ=s_r(Q)g_r+(q_r-q_r^{-1}){1\over e^r-1}(s_r(Q)-Q) }
where $Q$ is in the group algebra of the weight lattice.
So the above representation is simply obtained by considering the
action of $g_r$ on the group algebra $P$
which satisfy the affine
Hecke relations with $g_r$ and setting it equal to $q_r$ to the right
of an expression.

One important consequence of \repgg~ is that any polynomial
in $P$ which is Weyl invariant ($s_r(Q)=Q$ for all $r\in\Pi$)
is in the center of the affine Hecke algebra:
$[g_r,Q]=0$ for all $r\in\Pi$.

Finally, in order to obtain a representation of the operators $S_{\xi}$,
we must give  the realization of the operators $t^a$
which satisfies the relations \nyta.
The operators $t^a$
obey the following commutation
relations with $P$:
\eqn\fcomm{
t^ae^{\la}=e^{(\la,a)}e^{\la}t^a }

{\it remarks:}

In \repg , $g_r$ is expressed in terms of the reflection $s_r$,
alternatively, we can express the reflections $s_r$ in terms of the
generators $g_r$. We obtain an expression very similar to \reyy:
\eqn\reyyb{
s_r={e^r g_r-g_r^{-1}\over
e^r q_r-q_r^{-1}}}
but now, the Weyl group relations \uni,\unii~
satisfied by $s_r$  rely on the fact that
$e^r$ and $g_r$ obey the defining relations of the affine Hecke algebra.

In this representation, there are two sets of operators which obey
the affine Hecke relations with the operators $g_r$.
One is the group algebra $P$ generated by the spectral parameters
$e^{\xi}$ and the other is the group algebra generated by the $S_{\xi}$
computed as in \wordxi.

\newsec{Affine Hecke algebra, quadratic algebras and physical models}

Our aim is to relate the Affine-Hecke relations
to the quadratic algebras.
In this section we are concerned with specific
representations of the  quadratic algebras which use the
$x_r$. We show
that the spectral parameters
which enter the definition of $x_r$ in \rey~
can be taken to obey the affine-Hecke
relations \heka.
By considering the operators which commute with the
quadratic algebras, one
obtains a commuting set of Hamiltonians which
describe physical models.

In the first part, we give a brief description of
the quadratic algebras which we need in the following.
In the second part we describe the procedure of quantization
of the spectral parameters of the monodromy matrices which obey
the quadratic relations.
Finally, in the fourth part, we give a brief description of
the physical models which result from this construction.

\subsec{Quadratic algebras}

Let us consider the quadratic relations for the matrices $T_r$:
\eqn\bell{
x_{e_a-e_b}T_{e_a}T_{e_b}=
T_{e_b}T_{e_a}x_{e_a-e_b} }

and
\eqn\bella{
x_{e_a-e_b}T_{e_a}x_{e_a+e_b}T_{e_b}=
T_{e_b}x_{e_a+e_b}T_{e_a}x_{e_a-e_b} }

Note that these
quadratic algebras are obtained from the Yang-Baxter equations \ybb~
with $m_{r,s}=3,4$
by replacing  two of the generators by operators called "Monodromy matrices''.
The monodromy matrix $T_{e_a}$ can be expanded in the parameter $e^a$
called its specrtal parameter.
A realization of the relations \bell~ \bella~is respectively given by:
\eqn\reali{
T_{e_a}=x_{e_a-e_n}x_{e_a-e_{n-1}}...x_{e_a-e_1}  }
\eqn\real{
T_{e_a}=x_{e_a-e_n}x_{e_a-e_{n-1}}...
x_{e_a-e_1}x_{e_a}x_{e_a+e_1}...x_{e_a+e_n} }
The vectors $e_a,e_b,e_1,...,e_n$
form an orthogonal basis. We denote $V$ the vector space
spanned by the $e_i$, $i=1,...,n$.
$e_a,e_b$ are called  auxiliary vectors.
The relations \bell,\bella~result from the Yang-Baxter equations \ybb~
satisfied by the $x_r$.
In what follows, we consider the representation \rey~
of the operators $x_r$ in terms of generators $f_r$
of the Hecke algebra and $w_r$ of the Weyl group.
\eqn\reya{
x_r=w_r{e^{r}f_r-f_r^{-1}\over
e^{r}q_r-q_r^{-1}}  }
We recall that in this formula, the reflection $w_r$ acts on
the Hecke generator $f_r$ and that the spectral parameter
$e^r$ commutes with $f_r$ and $w_r$.

\subsec{Quadratic algebras  and affine Hecke algebras }

The Weyl group of $A_n$ ($S_n$),
acts by permuting the vectors  $e_i$ of $V$ in the
first case \bell~ and  The Weyl group of $B_n,C_n$ acts by permuting
the vectors and taking their opposite in the second case \bella.
In each case, we consider the fundamental systems
given by:
\eqn\fusi{
\Pi=\{e_1-e_2,...,e_{n-1}-e_n\} }
\eqn\fusil
{\Pi=\{e_1-e_2,...,e_{n-1}-e_n,e_n\} }

The generators $f_r$ with $r\in \Pi$
generate a Hecke algebra.
Let us consider another realization of the
Hecke algebra \hecb~ generated by $g_r$
with $r\in \Pi$. The $g_r$ commute with the $w_r$ and the $f_r$.

We quotient the group algebra of the $g_r$ and the $f_r$
by the relation $f_r=g_r$ to the right. We denote $\pi(.)$
the operation which consists in eliminating
the generators $g_r$ in the quotient: $\pi(...g_r)=\pi(...)f_r$.
We now require that the quadratic
relations \reali,\real, are still satisfied when one replaces
the monodromy matrices $T_{e_a}$ by their projection
$\pi(T_{e_a})$.
Let us  show that the relations are
satisfied if the spectral parameters $e^r$ obey
the defining relations \heka~ of the
affine Hecke algebra
with the generators $g_r$.

In order for the quadratic relations \bell~ \bella~ to be satisfied
after taking the quotient
The condition to satisfy is:
\eqn\gen{
\pi(g_rT_a)=f_r\pi(T_a) \ \ {\rm for}\ r\in\Pi}
It ensures that $\pi(T_aT_b)=\pi(T_a)\pi(T_b)$ and therefore that
the relations \bell~\bella~are satified
when one substitutes $\pi(T_a)$ for $T_a$ in them.

To verify \gen one commutes the generators
$g_r,\ r\in\Pi$ through $T_a$ using the
affine Hecke relations with the spectral
parameters which enters the definition  of the $x_r$.
Once $g_r$ has been pushed to the right,
one replaces it with $f_r$ and one commutes $f_r$ back to the left
using the Hecke relations and the fact that that the Weyl reflections
$w_r$ act on the $f_r$.
The computation is simplified by the fact that the denominator
of $T_a$ is a Weyl invariant polynomial in $P$. It therefore commutes with
the $g_r$ and one can keep
only the numerator of the $x_r$, $(x_r=w_r(e^rf_r-f_r^{-1}))$.

\subsec{Physical models}

To simplify the discussion, we restrict ourselves to the $A_{n-1}$ case \bell.
We consider $n$ particles with coordinates $\theta_j$ on the unit circle.
Each particle carries a spin denoted $\sigma_j$.
They are described by a wave function $\Psi(z_j)$ which is
supposed to be polynomial in the variables $z_j=e^{i\theta_j}$.

There is an action of the Hecke algebra on the coordinates which follows
from \repg:
\eqn\hou{
g_{j,j+1} =q s_{j,j+1} +(q-q^{-1}){z_{j+1}\over z_j -z_{j+1}}
(s_{j,j+1} -1) }
where the permutation $s_{j,k}$ permutes the coordinates $z_j$ and $z_k$.
There is also  a representation of the Hecke algebra on the spin variables
denoted $f_{j,j+1}$
($f_{j,j+1}$ is supposed to act on the spin of the particles $j$ and $j+1$).

The operators $t_i=t^{e_i}$ act by shifting the variable $\theta_i$:
\eqn\tpsi{
t_i\Psi(z_1,...,z_i,...)=\Psi(z_1,...,tz_i,...)  }

We call "Physical states" those for which the two Hecke algebras
act in the same way.
\eqn\phys{
g_{j,j+1}\Psi_{\rm phys}=
f_{j,j+1}\Psi_{\rm phys}  }
It is clear that on these states, we can replace the action of the operators
$g_{j,j+1}$
by $f_{j,j+1}$.
The operation which consists in  replacing the generators $g_{j,j+1}$
by $f_{j,j+1}$ coincides with the projection
$\pi(...g_{j,j+1})=\pi(...)f_{j,j+1}$
described before.
An operator $O$ acting on physical states preserves this space if:
$\pi(g_{j,j+1}O)=f_{j,j+1}O$.

Let us denote by $S_j$ the Hecke generators $S_{e_j}$ constructed in
\exaa. $S_j$ and $g_{k,k+1}$ obey the relations:
\eqn\hoka{
\eqalign{
[g_{k,k+1},S_j]&=0 \ {\rm if}\ j\ne k,k+1 \cr
g_{j,j+1}S_j&=S_{j+1}g_{j,j+1}^{-1}  \cr }}

We consider the following Hamiltoniens acting on physical states:
\eqn\hamol{
H_l=\sum_{1\le i_1<i_2..<i_l\le n} S_{i_1}...S_{i_l}  }

These operators obviously commute with each other. Moreover,
since they are symmetric polynomials in the $S_j$,
they
commute with $g_{j,j+1}\ \forall j$.
Thus, their action preserves the physical states.
In the simple case where $f_{j,j+1}=q$ (or $-q^{-1}$)
(the scalar case), the projection of
these Hamiltoniens by $\pi$ coincides with the trigonometric models
defined by Ruijenhaars.
In this case, the Hamiltonians are invariant under the permutations
of the coordinates which makes the identification with Ruijenhaars
models easy to do.
Unfortunately, in the general case,
this projection is much more difficult to describe explicitly.

One can also construct the monodromy matrix \reali.
The expression of $x_{e_a-e_i}$ being:
\eqn\nexi{
x_{e_a-e_i}=s_{a,i}{e^uf_{e_a-e_i}-S_if_{e_a-e_i}^{-1}
\over e^uq-S_iq^{-1}} }
Here $e^u$ is the spectral parameter of the monodromy matrix.
Because of \gen,
the monodromy matrix
also preserves the physical states and it commutes with $H_l$.
This remains true when one replaces the operators by their
projection by $\pi$.
Thus, the algebra defined by the quadratic relations \bell~ is a
symmetry algebra for the $H_l$.

Let us show that the operators $S_j$ constructed in this way
can be diagonalized simultaneously and are unitary.

{\it Eigenvalues of the } $S_j$:

We show that the operators $S_j$ are represented by triangular matrices.
Let us recall the expression for the $S_j$:
\eqn\ess{
S_j=x_{j-1,j}^{-1}...
x_{1,j}^{-1} t_j
x_{j,n}...x_{j,j+1}  }

where the operator $x_{i,j}$ takes the limiting form for $i<j$:
\eqn\limf{
x_{i,j}=q+(q-q^{-1}){z_i \over (z_i-z_j)} (s_{i,j}-1) }
$x_{i,j}$ commutes with $z_iz_j$ and with $z_k$ for $k\ne i,j$.
It acts in the following way for on the monomials $z_i^m,\ z_j^m
,\ m>0$ and $1$:
\eqn\act{
\eqalign{
x_{i,j}z_i^m&=q^{-1}z_i^m-(q-q^1)(z_i^{m-1}z_j+...+z_iz_j^{m-1})\
{\rm for}\  m>0\cr
x_{i,j}z_j^m&=q z_j^m+(q-q^{-1})(z_i^m+...+z_iz_j^{m-1})\  {\rm for}\  m>0\cr
x_{i,j}z_i^0&=0\cr }}
To a monomial $z_1^{k_1}...z_n^{k_n}$
we associate a partition $|k|=(k_{p_1}\ge k_{p_2}...\ge k_{p_n})$
where we order the $k_j$ in decreasing order.

Let us consider which
new monomials
$z_1^{k_1'}...z_n^{k_n'}$ can
appear
when one acts $x_{i,j}$ on this monomial.
First, all the $k'_l$ for $l\ne i,j$ are equal to $k_l$.
Then, if the partition of the new monomial
is different from $|k|$, it can be obtained
from $|k|$ by a sequence of squeezing operations:
$(...,k_i,...,k_j,....)\to (...,k_i-1,...,k_j+1,...)$  if $k_i>k_j$.
$(...,k_j,...,k_i,....)\to (...,k_j-1,...,k_i+1,...)$  if $k_i<k_j$.
Finally if the partition of the new monomial is equal to $|k|$,
$k'_i=k_i$ and $k'_j=k_j$ if $k_i>k_j$.

Let us  define an order on the  monomials by saying that $z_i^{k_i}$
is larger than $z_i^{k'_i}$ if either $|k'|$ is obtained from
$|k|$ by a sequence of squeezing operations, or $|k|=|k'|$
and $k'_i-k_i>0\ \forall i$.
It follows from the above analyses that the action of $S_j$ on a monomial
produces only monomials which are smaller with respect to this order.
Thus the eigenvalues of the operators $S_j$ are given by the diagonal elements
in the monomial basis.

Given the partition $|k|=(k_1,...,k_n)$, the eigenvalues corresponding to the
monomials associated to it are  all obtained by
permutations of the multiplet:
\eqn\multiplet{
(S_j)=(t^{k_j}q^{n+1-2j}) }
If we set $q=t^{\beta/2}$ with $\beta$ a real parameter, the
operators $S_j$ have the physical interpretation of exponentials
of momentum operators: $S_j=t^{K_j}$.
$K_j=k_j+\beta({n+1\over2}-j)$ obey a generalized Pauli principle
since they must be $\beta$ apart from each other.

{\it Scalar product}

A scalar product can be defined \cher ~, so that the operators
$g_{j,j+1}$ and $S_j$ are  unitary:
$g_{j,j+1}^+=g_{j,j+1}^{-1}$,
$S_j^+=S_j^{-1}$.
For $q=t^{k/2}$, with $k$ an integer,
the scalar product is given by:
\eqn\prosca{
<\Psi_1,\Psi_2>=\int \prod_{i=1}^n d\theta _i \ \overline\Psi_1(z_j)
\Psi_2(z_j) C(z_j) }
where the bar symbol stands for the complex conjugation.
$q,t,z_j$ are supposed to be
complex numbers of modulus 1 and
the integration over $\theta_j$ keeps the coefficient of $z^0$ of
the integrand.

The measure $C(z_k)$ is given by:
\eqn\mesure{
C(z_j)=\prod_{i<j} \prod_{l=-k}^{k-1} \ (\sqrt{t^lz_i/z_j}
-\sqrt{t^{-l}z_j/z_i}) }

In the case of the Ruijenhaars models, many properties of the wave
functions $\Psi$ can be obtained, for example their norms.
We refer to \mac~ for a more complete analyses of their properties.

\listrefs
\end